\documentclass{amsart}
\usepackage{amssymb}
\usepackage{graphics}
\usepackage{amscd}
\newtheorem{thm}{Theorem}[section]
\newtheorem{cor}[thm]{Corollary}
\newtheorem{lem}[thm]{Lemma}
\newtheorem{prop}[thm]{Proposition}
\theoremstyle{definition}
\newtheorem{defn}[thm]{Definition}

\theoremstyle{remark}
\newtheorem{rem}[thm]{Remark}
\numberwithin{equation}{section}
\newtheorem{note}[thm]{Notation}
\newcommand{\abs}[1]{\left\vert#1\right\vert}

\newcommand{\Real}{\mathbb R}
\newcommand{\eps}{\varepsilon}

\newcommand{\im}{\mathop{\textrm{Im}}\nolimits}
\newcommand{\imm}{\mathop{\textrm{Imm}}\nolimits}
\begin{document}

\title{A geometric classification of immersions of 3-manifolds into 5-space}%
\author{Andr\'as Juh\'asz}%
\address{Department of Analysis, E\"otv\"os Lor\'and University,
P\'azm\'any P\'eter s\'et\'any 1/C, Budapest, Hungary 1117}%
\email{juhasz.6@dpg.hu}%

% \thanks{Research partially supported by OTKA grant no. T037735}%

\date{\today}%
%\dedicatory{}%
%\commby{}%
% ----------------------------------------------------------------
\begin{abstract}
In this paper we define two regular homotopy invariants $c$ and
$i$ for immersions of oriented 3-manifolds into $\Real^5$ in a
geometric manner. The pair $(c(f),i(f))$ completely describes the
regular homotopy class of the immersion $f$. The invariant $i$
corresponds to the 3-dimensional obstruction that arises from
Hirsch-Smale theory and extends the one defined in \cite{Takase}
for immersions with trivial normal bundle.
\end{abstract}
\maketitle
% ----------------------------------------------------------------
\section{Introduction}

Hirsch-Smale theory \cite{Hirsch} reduces the problem of regular
homotopy classification of immersions to homotopy theory. However,
this homotopy theoretic problem is usually hard to deal with. In
the case of immersions of oriented 3-manifolds into $\Real^5$ this
homotopy theoretic problem was solved by Wu \cite{Wu} using
algebraic topological methods (also see \cite{Li}).

However, it remains a problem to determine the regular homotopy
class of a given immersion from its geometry. By geometry we mean
the structure of the "singularities" of the map. For example,
double points are such "singularities", and indeed, Smale
\cite{Smale2} showed that for $n > 1$ the regular homotopy class
of an immersion $S^n \looparrowright \Real^{2n}$ is completely
determined by the number of its double points (modulo 2 if $n$ is
odd). A similar classification was carried out by Ekholm
\cite{Ekholm2} for immersions of $S^k$ into $\Real^{2k-1}$ for $k
\ge 4$.

The whole picture changes when we consider immersions of $S^3$
into $\Real^5$. Hughes and Melvin \cite{Hughes} showed that there
are infinitely many embeddings $S^3\hookrightarrow \Real^5$ that
are pairwise not regularly homotopic to each other. Therefore one
can not determine the regular homotopy class from the
"singularities" since an embedding has no such. Ekholm and
Sz\H{u}cs \cite{ESz} came over this problem using "singular
Seifert surfaces" bounded by the immersions. For an immersion $f
\colon M^3 \looparrowright \Real^5$ a singular Seifert surface is
a generic map $F \colon W^4 \to \Real^5$ of a compact orientable
manifold $W^4$ with boundary $M^3$ such that $\partial F = f$. In
\cite{ESz} it is shown that for $M^3 = S^3$ the Smale invariant of
$f$ can be computed from the singularities of $F$. Later Saeki,
Sz\H{u}cs and Takase \cite{Takase} generalized these results for
immersions $f \colon M^3 \looparrowright \Real^5$ with trivial
normal bundle (for oriented $M^3$). The invariant introduced in
\cite{Takase} corresponds to the 3-dimensional obstruction to a
regular homotopy between two such immersions. Our present paper
generalizes the results of \cite{Takase} to arbitrary immersions
$f \colon M^3 \looparrowright \Real^5$.

We will consider the set $\imm(M^3, \Real^5)_{\chi}$ of immersions
with fixed normal Euler class $e(\nu_f) = \chi \in H^2(M^3;
\mathbb{Z})$ and construct a $\mathbb{Z}_{2d(\chi)}$-valued
regular homotopy invariant $i$ for this set of immersions, where
$d(\chi)$ denotes the divisibility of $\chi$. The construction of
the invariant $i$ will also make use of a singular Seifert surface
$F$. In \cite{Takase} $F$ had to be an immersion near the
boundary, but we (have to and) will allow arbitrary generic maps.
If $\chi = 0$ and $F$ is an immersion near the boundary then the
construction of the invariant $i$ agrees with the one introduced
in \cite{Takase}. We will also show that whenever $f, g \colon M^3
\looparrowright \Real^5$ are regularly homotopic on a neighborhood
of the 2-skeleton of $M^3$ then $i(f) = i(g)$ iff $f$ and $g$ are
regularly homotopic. This shows that $i$ corresponds to the
3-dimensional obstruction to a regular homotopy between $f$ and
$g$. (Note that there is an invariant which determines the regular
homotopy class of the restriction of an immersion to a
neighborhood of the 2-skeleton of $M^3$. This invariant was called
the Wu invariant in \cite{Takase}, see below.)

Regular homotopy classes of immersions of oriented 3-manifolds
into $\Real^5$ endowed with the connected sum operation form a
semigroup whose structure we will also determine. Finally, an
exact sequence will be defined that relates $\imm[M^3, \Real^5]$
to $\imm[M^3, \Real^6]$ and $[M^3, S^2]$.

\section{Preliminaries}

First let us recall the result of Wu \cite{Wu} that classifies
immersions of an oriented 3-manifold $M^3$ into $\Real^5$ up to
regular homotopy.

\begin{thm} \label{thm:1}
The normal Euler class $\chi$ of an immersion $f \colon M^3
\looparrowright \Real^5$ is of the form $2c$ for some $c \in
H^2(M^3; \mathbb{Z})$ and for any $c \in H^2(M^3; \mathbb{Z})$
there is an immersion $f$ such that $\chi = 2c$. Furthermore,
$$\imm[M^3,\Real^5]_{\chi} \approx \coprod_{c \in H^2(M^3;
\mathbb{Z}) \,,\, 2c = \chi} H^3(M^3; \mathbb{Z})/(2\chi \cup
H^1(M^3; \mathbb{Z})),$$ where $\imm[M^3,\Real^5]_{\chi}$ is the
set of regular homotopy classes of immersions with normal Euler
class $\chi \in H^2(M^3; \mathbb{Z})$ and $\cup$ represents the
cup product, moreover the symbol $\approx$ denotes a bijection.
\end{thm}

\begin{rem}
For $\chi \in H^2(M^3; \mathbb{Z})$ let $d(\chi) \in \mathbb{Z}$
denote the divisibility of $\chi$, so that $\chi$ equals $d(\chi)$
times a primitive class in $H^2(M^3; \mathbb{Z})$ modulo torsion,
and $d(\chi) = 0$ if $\chi$ is of finite order. Then Poincar\'e
duality implies that $$H^3(M^3; \mathbb{Z})/(2\chi \cup H^1(M^3;
\mathbb{Z})) \approx \mathbb{Z}_{2d(\chi)}.$$ If $f$ is an
immersion of $M^3$ into $\Real^5$ with normal Euler class $\chi$
then let us introduce the notation $d(f)$ for $d(\chi)$.
\end{rem}

\begin{note}
For $\chi \in H^2(M^3; \mathbb{Z})$ let $\Gamma_2(\chi)$ denote
the set $\{\,c \in H^2(M^3; \mathbb{Z}) \colon 2c = \chi \,\}$.
Throughout this paper we will use the notation $M^3_{\circ}$ for
the punctured 3-manifold $M^3 \setminus D^3$, where $D^3 \subset
M^3$ is a closed 3-disc. Then the 2-skeleton $\text{sk}_2(M^3)$ is
a deformation retract of $M^3_{\circ}$.
\end{note}

Theorem \ref{thm:1} can also be applied to the open manifold
$M^3_{\circ}$. Since $H^3(M^3_{\circ}; \mathbb{Z}) = 0$ we obtain
a bijection $$\bar{c} \colon \imm[M^3_{\circ}, \Real^5]_{\chi} \to
\Gamma_2(\chi).$$ Thus for an immersion $f \colon M^3
\looparrowright \Real^5$ the invariant $c(f) =
\bar{c}(f|M^3_{\circ}) \in \Gamma_2(\chi)$ describes the regular
homotopy class of $f|M^3_{\circ}$. Following \cite{Takase} we will
call $c(f)$ the Wu invariant of the immersion $f$.

To get a complete description of $\imm[M^3, \Real^5]_{\chi}$ we
will construct a $\mathbb{Z}_{2d(\chi)}$-valued invariant $i$ such
that the map $$(c, i) \colon \imm[M^3, \Real^5]_{\chi} \to
\Gamma_2(\chi) \times \mathbb{Z}_{2d(\chi)}$$ will be a bijection.
The invariant $i$ is constructed in a geometric manner and is an
extension of the invariant defined in \cite{Takase} for $\chi =
0$.

Next let us recall Theorem 1.1(a) in \cite{ESz}. Let $f \colon S^3
\looparrowright \Real^5$ be an immersion and $V^4$ an arbitrary
compact oriented 4-manifold with $\partial V^4 = S^3$. The map $f$
extends to a generic map $F \colon V^4 \to \Real^5$ which has no
singular points near the boundary $\partial V^4$ since the normal
bundle $\nu_f$ of $f$ is trivial. This map $F$ has isolated cusps,
each one having a sign. Let us denote by $\#\Sigma^{1,1}(F)$ their
algebraic number and let $\Omega(f)$ be the Smale invariant of
$f$. The following formula was proved in \cite{ESz}.

\begin{thm} \label{thm:2}
$$\Omega(f) = \frac12(3\sigma(V^4) + \#\Sigma^{1,1}(F)).$$
\end{thm}

The proof of this theorem relies on the following proposition
(\cite{Szucs}, Lemma 3).

\begin{lem} \label{lem:0}
Let $X^4$ be a closed oriented 4-manifold and $g \colon X^4 \to
\Real^5$ a generic map. Then $3 \sigma(X^4) + \#\Sigma^{1,1}(g) =
0$.
\end{lem}

For the sake of completeness we will also recall from
\cite{Takase} the definition of the invariant $i$ for immersions
with trivial normal bundle. First we need a preliminary
definition.

\begin{defn}
Let $M^3$ be a closed oriented 3-manifold. We denote by
$\alpha(M^3)$ the dimension of the $\mathbb{Z}_2$ vector space
$\tau H_1(M^3; \mathbb{Z}) \otimes \mathbb{Z}_2$, where $\tau
H_1(M^3; \mathbb{Z})$ is the torsion subgroup of $H_1(M^3;
\mathbb{Z})$.
\end{defn}

\begin{defn} \label{defn:6}
Let $f \colon M^3 \looparrowright \Real^5$ be an immersion with
trivial normal bundle. Let $W^4$ be any compact oriented
4-manifold with $\partial W^4 = M^3$ and $F \colon W^4 \to
\Real^5$ a generic map nonsingular near the boundary such that
$F|\partial W^4 = f$. (We can choose such a generic map $F$ since
$f$ is an immersion with trivial normal bundle.) Denote the
algebraic number of cusps of $F$ by $\#\Sigma^{1,1}(F)$. Then let
$$i(f) = \frac32 (\sigma(W^4)-\alpha(M^3)) + \frac12\#\Sigma^{1,1}(F).$$
It is proved in \cite{Takase} that $i(f)$ is always an integer and
a regular homotopy invariant.
\end{defn}

In the following sections we will extend the above regular
homotopy invariant $i$ to arbitrary immersions. If $f \colon M^3
\looparrowright \Real^5$ has non-trivial normal bundle then we
have to give up the assumption that the singular Seifert-surface
$F$ is an immersion near the boundary. Thus we will use an
arbitrary generic map $F \colon W^4 \to \Real^5$ such that
$\partial F = f$. The singular set $\Sigma^1(F)$ of such an $F$ is
a 2-dimensional submanifold of $W^4$ with boundary $C(F) =
\partial \Sigma^1(F) \subset M^3$. If we orient $\ker(dF)|C(F)$
so that it points into $W^4$ and project it into $TM^3$ then we
obtain a normal field $\nu(F)$ along $C(F)$. We will define the
rotation of $\nu(F)$ around $C(F)$ modulo $4d(f)$ and denote this
by $R(F)$. The double of the extended invariant will be defined to
be
$$I(f) = 3(\sigma(W^4) - \alpha(M^3)) + \#\Sigma^{1,1}(F) + R(F)
\in \mathbb{Z}_{4d(f)}.$$ We will show that $I(f)$ is always even,
thus it defines an element $i(f) \in \mathbb{Z}_{2d(f)}$ using the
natural embedding $\mathbb{Z}_{2d(f)} \hookrightarrow
\mathbb{Z}_{4d(f)}$.

\section{Rotation}

Throughout this paper $M^3$ will denote a fixed closed connected
and oriented 3-manifold.

\begin{note} \label{note:1}
A pair $(C, \nu)$ will always stand for an oriented 1-dimensional
submanifold $C$ of $M^3$ and a nowhere vanishing normal field
$\nu$ along $C$.
\end{note}

\begin{defn} \label{defn:1}
Let $\chi \in H^2(M^3; \mathbb{Z})$ and let $C_0$ and $C_1$ be
1-dimensional oriented submanifolds of $M^3$ with normal fields
$\nu_0$ and $\nu_1$ such that $PD[C_0] = PD[C_1] = \chi$. (Here
$PD$ denotes Poincar\'e duality.) Then we can define the
\emph{rotation difference} $\text{rd}((C_0,\nu_0),(C_1,\nu_1)) \in
\mathbb{Z}_{2d(\chi)}$ of $(C_0, \nu_0)$ and $(C_1, \nu_1)$ as
follows. Since $[C_0] = [C_1]$ and
$$H_1(M^3; \mathbb{Z}) \approx H^2(M^3; \mathbb{Z}) \approx [M^3,
\mathbb{C}P^{\infty}],$$ there exists an oriented cobordism $K^2
\subset M^3 \times I$ between $C_0 \subset M^3 \times\{0\}$ and
$C_1 \subset M^3 \times \{1\}$. Let $\nu$ be a generic normal
field along $K^2$ that extends $\nu_0$ and $\nu_1$. Then a sign
can be given to each zero of $\nu$ since $M^3$ is oriented. Now we
define $\text{rd}((C_0,\nu_0),(C_1,\nu_1))$ to be the algebraic
number of zeroes of $\nu$ modulo $2d(\chi)$. Equivalently,
$\text{rd}((C_0,\nu_0),(C_1,\nu_1))$ is the self intersection of
$K$ in $M^3 \times I$ modulo $2d(\chi)$ if perturbed in the
direction of $\nu$.
\end{defn}

\begin{rem} \label{rem:3}
The rotation difference is the obstruction to the existence of a
framed cobordism between the framed submanifolds $(C_0, \nu_0)$
and $(C_1, \nu_1)$ of $M^3$. Using the Pontrjagin construction
this corresponds to the obstruction to a homotopy between two maps
of $M^3$ to $S^2$. This situation was first examined in
\cite{Pontrjagin}. It is easy to see that $\text{rd}((C_0,
\nu_0),(C_1, \nu_1)) = 0$ iff $(C_0, \nu_0)$ and $(C_1, \nu_1)$
are framed cobordant. Thus we obtain a bijection $$[M^3, S^2]
\approx \coprod_{\chi \in H^2(M^3; \mathbb{Z})} H^3(M^3;
\mathbb{Z})/ 2\chi \cup H^1(M^3 ; \mathbb{Z}).$$
\end{rem}

\begin{prop} \label{prop:1}
In Definition \ref{defn:1} above the rotation difference is well
defined, i.e., it does not depend on the choice of $K$ and $\nu$.
\end{prop}

\begin{proof}
Let $K$, $\nu$ and $K'$, $\nu'$ be as in Definition \ref{defn:1}.
We glue together $M^3 \times I$ and $-M^3 \times I$ along their
boundaries so that we obtain the double $D(M^3 \times I) = M^3
\times S^1$. Place $K$ into the half of $M^3 \times S^1$
corresponding to $M^3 \times I$ and $K'$ into the other half. Then
we obtain a closed oriented surface $F = K \cup -K'$ in $M^3
\times S^1$ and a normal field $\mu = \nu \cup \nu'$ along $F$.
Since $H^*(S^1; \mathbb{Z})$ is a torsion free $\mathbb{Z}$-module
we can apply K\"unneth's theorem and we get that
$$H^2(M^3 \times S^1; \mathbb{Z}) \approx H^1(M^3; \mathbb{Z}) \otimes
H^1(S^1; \mathbb{Z}) \oplus H^2(M^3; \mathbb{Z})\otimes H^0(S^1;
\mathbb{Z}).$$  Thus the Poincar\'e dual of $F$ can be written in
the form $$ PD[F] = x \times \alpha + y \times 1 \in H^2(M^3
\times S^1; \mathbb{Z}),$$ where $x \in H^1(M^3; \mathbb{Z})$ and
$y \in H^2(M^3; \mathbb{Z})$, moreover $\alpha$ denotes the
generator of $H^1(S^1; \mathbb{Z})$ and $1$ the generator of
$H^0(S^1; \mathbb{Z})$ given by the orientation of $S^1$. Note
that for $1 \in S^1$ the dual class of $M^3 \times \{1\} \subset
M^3 \times S^1$ is $$PD[M^3 \times \{1\}] = PD[M^3] \times
PD[\{1\}] = 1 \times \alpha \in H^1(M^3 \times S^1; \mathbb{Z}).$$
Moreover,
$$PD[F \cap (M^3 \times \{1\})] = PD[C_1 \times \{1\}] = PD[C_1]
\times PD[\{1\}] = \chi \times \alpha.$$ On the other hand
$$PD[F \cap (M^3 \times \{1\})] = PD[F] \cup PD[M^3 \times \{1\}] =
(x \times \alpha + y \times 1) \cup (1 \times \alpha) = x \times
\alpha^2 + y \times \alpha.$$ Since $\alpha^2 = 0$ we get that $y
\times \alpha = \chi \times \alpha$. Using K\"unneth's theorem
again we obtain the equality $y = \chi$. Thus we get that $$PD[F]
\cup PD[F] = (x \times \alpha + \chi \times 1)^2 = (2x \cup \chi)
\times \alpha $$ since $\alpha^2 = \chi^2 = 0$ and $x \cup \chi =
\chi \cup x$ because the degree of $\chi$ is 2. So the self
intersection of $F$ in $M^3 \times S^1$ equals $\langle (2\chi
\cup x) \times \alpha, [M^3 \times S^1] \rangle = \langle 2\chi
\cup x, [M^3] \rangle \in 2d(\chi)\mathbb{Z}$. If we perturb $F$
in the direction of $\mu$ we get that the self intersection of $K$
with respect to $\nu$ equals the self intersection of $K'$ with
respect to $\nu'$ modulo $2d(\chi)$.
\end{proof}

\begin{prop}
If $[C_0]=[C_1]=[C_2] \in H_1(M^3, \mathbb{Z})$ then
$$\text{rd}((C_0,\nu_0),(C_1,\nu_1)) + \text{rd}((C_1,\nu_1),(C_2,\nu_2)) =
\text{rd}((C_0,\nu_0),(C_2,\nu_2)).$$
\end{prop}

\begin{defn}
For each $a \in H_1(M^3; \mathbb{Z})$ fix a pair $(C_a, \nu_a)$
such that $[C_a] = a$. Then for $[C] = a$ let $r(C,\nu) =
\text{rd}((C,\nu),(C_a, \nu_a))$.
\end{defn}

\begin{cor}
If $[C_0]=[C_1]$ then $r(C_0,\nu_0) - r(C_1, \nu_1) =
\text{rd}((C_0, \nu_0), (C_1, \nu_1))$.
\end{cor}

\begin{defn} \label{defn:7}
We can define the mod 2 rotation difference $\text{rd}_2((C_0,
\nu_0), (C_1, \nu_1))$ for unoriented $C_0$ and $C_1$ just as in
Definition \ref{defn:1} but allowing the cobordism $K$ to be
non-orientable and counting the self intersection of $K$ in $M
\times I$ only modulo 2. The proof that this is well defined is
analogous to the oriented case. It is clear that the epimorphism
$\mathbb{Z}_{2d(\chi)} \to \mathbb{Z}_2$ takes $\text{rd}$ to
$\text{rd}_2$. The mod 2 rotation $r_2$ is defined just like $r$.
\end{defn}

Unfortunately we will have to lift the invariants $\text{rd}$ and
$r$ to $\mathbb{Z}_{4d(\chi)}$. To be able to do this we need more
structure on $M^3$ then just a framed submanifold. We will use
this additional structure to restrict the homology class of the
cobordism $K$ so that the surface $F$ in the proof of Proposition
\ref{prop:1} will represent an even homology class and thus $x$
will always be even (since $\chi$ is even). So the self
intersection of $F$ will be divisible by $4d(\chi)$ instead of
just $2d(\chi)$.

\begin{note}
Fix a cohomology class $\chi \in H^2(M^3; \mathbb{Z})$. Let
$\eps^3_M$ denote the 3-dimensional trivial bundle over $M^3$ and
let $t,v \in \Gamma(\eps^3_M)$ be two generic non-zero sections of
$\eps^3_M$. Furthermore, suppose that the 2-dimensional oriented
subbundle $t^{\perp} < \eps^3_M$ has Euler class $\chi$. If we
project $v$ into $t^{\perp}$ we obtain a section $w \in
\Gamma(t^{\perp})$ that vanishes along a curve $C \subset M^3$ and
we orient $C$ so that $PD[C] = e(t^{\perp})$. In particular, $t$
and $v$ are linearly dependent exactly at the points of $C$.
Finally let $\nu$ be a non-zero normal field along $C$. In the
future we will denote such a structure on $M^3$ by a quadruple
$(C, \nu, t,v)$ and the set of these structures by $N(M^3, \chi)$.
\end{note}

\begin{rem}
Since $PD[C]|_2 = w_2(\eps^3_M) = 0 \in H^2(M^3; \mathbb{Z}_2)$,
the cohomology class $\chi = PD[C]$ is of the form $2c$ for some
$c \in H^2(M^3; \mathbb{Z})$. This can be seen from the long exact
sequence associated to the coefficient sequence $\mathbb{Z} \to
\mathbb{Z} \to \mathbb{Z}_2$. Thus $N(M^3, \chi) = \emptyset$ if
$\chi$ is not of the form $2c$.
\end{rem}

\begin{defn} \label{defn:5}
Suppose that $a_0 = (C_0, \nu_0, t_0, v_0)$ and $a_1 = (C_1,
\nu_1, t_1, v_1)$ are elements of $N(M^3, \chi)$, where $\chi =
2c$. Then we will define their rotation difference
$\text{Rd}(a_0,a_1) \in \mathbb{Z}_{4d(\chi)}$ as follows. We will
consider $a_i$ to be in $N(M^3 \times \{i\}, \chi)$ for $i = 0,1$.
Let $t, v \in \Gamma(\eps^3_{M \times I})$ be generic non-zero
sections extending $t_i$ and $v_i$ for $i = 0,1$. Denote by $K$
the 2-dimensional submanifold of $M^3 \times I$ where $t$ and $v$
are linearly dependent. Let $w$ denote the projection of $v$ into
the 2-dimensional oriented subbundle $t^{\perp} < \eps^3_{M \times
I}$. Then $w$ is zero exactly at the points of $K$, thus it
defines an orientation of $K$. With this orientation $K$ is an
oriented cobordism between $C_0$ and $C_1$. Let $\nu$ denote a
normal field of $K$ that extends both $\nu_0$ and $\nu_1$.  Now we
define $\text{Rd}(a_0,a_1)$ to be the algebraic number of zeroes
of $\nu$ modulo $4d(\chi)$. Equivalently, $\text{Rd}(a_0,a_1)$ is
the self intersection of $K$ in $M^3 \times I$ modulo $4d(\chi)$
if perturbed in the direction of $\nu$.
\end{defn}

\begin{prop} \label{prop:2}
In Definition \ref{defn:5} the rotation difference is well
defined. I.e., it does not depend on the extensions $t,v$ and
$\nu$.
\end{prop}

\begin{proof}
Let $t,v,\nu$ and $t',v',\nu'$ be as in Definition \ref{defn:5}.
The sections $t,v$ are linearly dependent over $K$ and $t',v'$ are
dependent over $K'$. Just as in the proof of Proposition
\ref{prop:1} we will place $K, \nu$ and $ K', \nu'$ in the two
halves of the double $D(M^3 \times I) = M^3 \times S^1$ and place
$t,v$ and $t',v'$ in the two halves of the trivial bundle
$\eps^3_{M \times S^1}$. Let $F$ denote the oriented surface $K
\cup -K'$ and by $\mu$ the normal field along $F$ obtained from
$\nu$ and $\nu'$. Moreover, let $T = t \cup t'$ and $V = v \cup
v'$. Then $T,V \in \Gamma(\eps^3_{M \times S^1})$ are linearly
dependent exactly over $F$, thus $PD[F]|_2 = w_2(\eps^3_{M \times
S^1}) = 0 \in H^2(M^3 \times S^1; \mathbb{Z}_2)$. Using the
coefficient sequence $\mathbb{Z} \to \mathbb{Z} \to \mathbb{Z}_2$
we get that $PD[F]$ is of the for $2b$ for some $b \in H^2(M^3
\times S^1; \mathbb{Z})$. Since $PD[F]$ is of the form $x \times
\alpha + \chi \times 1$ where $\chi = 2c$ we get that there exists
an element $z \in H^1(M^3; \mathbb{Z})$ such that $x = 2z$. Thus
$PD[F]^2 = (4z \cup \chi) \times \alpha$ which implies that the
self intersection of $F$ is divisible by $4d(\chi)$.  If we
perturb $F$ in the direction of $\mu$ we get that the self
intersection of $K$ with respect to $\nu$ equals the self
intersection of $K'$ with respect to $\nu'$ modulo $4d(\chi)$.
\end{proof}

\begin{rem} \label{rem:1}
The surface $K$ represents the dual of the Stiefel-Whitney class
of the bundle $\eps^3_{M \times I}$ relative to the sections $t_i,
v_i$ given over $M^3 \times \{0,1\}$. I.e., $$PD[K]|_2 =
w_2(\eps^3_{M \times I} ; t_i, v_i) \in H^2(M^3 \times I, ((M^3
\setminus C_0) \times \{0\}) \cup ((M^3 \setminus C_1)\times
\{1\}); \mathbb{Z}_2)$$ since $v$ and $t$ are linearly independent
over $((M^3 \setminus C_0) \times \{0\}) \cup ((M^3 \setminus
C_1)\times \{1\})$. Using Lefschetz duality we get that the
relative homology class $[K]|_2 \in H^2(M^3 \times I, C_0 \times
\{0\} \cup C_1 \times \{1\}; \mathbb{Z}_2)$ is independent of the
choice of $t$ and $v$. If we choose a simplicial subdivision of
$M^3$ so that $\text{sk}_1(M^3) \cap C_i = \emptyset$ for $i =
0,1$ then $w_2$ is the obstruction to extending the map $(t_i,v_i)
\colon \text{sk}_1(M^3 \times \{0,1\}) \to V_2(\Real^3)$ to
$\text{sk}_2(M^3 \times I)$. So the homology class $[K]|_2$ and
thus $\text{Rd}(a_0,a_1)$ depends only on the homotopy class of
the map
$$(t_i,v_i)|\text{sk}_1(M^3) \colon \text{sk}_1(M^3) \to V_2(\Real^3)$$ for $i =
0,1$. For the sake of completeness we note that if the extension
$t$ is given then
$$PD[K] = e(t^{\perp}; w_i) \in H^2(M^3 \times I, (M^3 \setminus
C_0) \times \{0\} \cup (M^3 \setminus C_1)\times \{1\};
\mathbb{Z}).$$ So we have obtained the following proposition.
\end{rem}

\begin{prop} \label{prop:3}
Suppose that $(C_0,\nu_0)$ and $(C_1, \nu_1)$ are framed
submanifolds of $M^3$ and let $a_0,b_0,a_1,b_1 \in N(M^3, \chi)$
be of the form $a_i = (C_i, \nu_i, t^a_i,v^a_i)$ and $b_i = (C_i,
\nu_i, t^b_i, v^b_i)$ for $i = 0,1$. Moreover, suppose that
$\text{sk}_1(M^3) \cap C_i = \emptyset$ and $(t^a_i,
v^a_i)|\text{sk}_1(M^3)$ is homotopic to $(t^b_i,
v^b_i)|\text{sk}_1(M^3)$ as maps into $V_2(\Real^3)$ for $i =
0,1$. Then the following equality holds:
$$\text{Rd}(a_0,a_1) = \text{Rd}(b_0,b_1).$$
\end{prop}

\begin{proof}
Let $t^a$ and $v^a$ be generic extensions of $t_i^a$, respectively
$v_i^a$ over $M^3 \times I$ and denote by $K^a$ the submanifold of
$M^3 \times I$ where $t^a$ and $v^a$ are linearly dependent. We
obtain the sections $t^b$ and $v^b$ of $\eps^3_{M \times I}$ and
the submanifold $K^b \subset M^3 \times I$ in a similar way. Then,
according to Remark \ref{rem:1}, we get that $$PD[K^a]|_2 =
w_2(\eps^3_{M \times I}; t^a_i,v^a_i) = w_2(\eps^3_{M \times I};
t^b_i,v^b_i) = PD[K^b]|_2,$$ since $w_2$ is the obstruction to
extending a map into $V_2(\Real^3)$ from $\text{sk}_1(M^3 \times
I)$ to $\text{sk}_2(M^3 \times I)$ and $(t^a_i,
v^a_i)|\text{sk}_1(M^3 \times \{i\})$ is homotopic to $(t^b_i,
v^b_i)|\text{sk}_1(M^3 \times \{i\})$. Thus if $F$ denotes the
submanifold of $M^3 \times S^1 = D(M^3 \times I)$ obtained by
piecing together $K^a$ and $K^b$ we get that $PD[F]|_2 =
w_2(\eps^3_{M \times S^1}) = 0$, so we can proceed as in the proof
of Proposition \ref{prop:2}.
\end{proof}

\begin{prop} \label{prop:5}
If $a_i = (C_i,\nu_i, t_i,v_i) \in N(M^3, \chi)$ for $i = 0,1$
then $$\text{rd}((C_0,\nu_0),(C_1,\nu_1)) \equiv
\text{Rd}(a_0,a_1) \mod 2d(\chi).$$
\end{prop}

\begin{prop}
If $a_0,a_1,a_2 \in N(M^3, \chi)$ then $$\text{Rd}(a_0,a_1) +
\text{Rd}(a_1,a_2) = \text{Rd}(a_0,a_2).$$
\end{prop}

\begin{defn}
For each $\chi \in H^2(M^3; \mathbb{Z})$ of the form $\chi = 2c$
fix an element $a_{\chi} \in N(M^3, \chi)$. Then for each $a \in
N(M^3, \chi)$ define the dotation $R(a) \in \mathbb{Z}_{4d(\chi)}$
to be $\text{Rd}(a,a_{\chi})$.
\end{defn}

\begin{cor}
If $a_0, a_1 \in N(M^3, \chi)$ then $\text{Rd}(a_0,a_1) = R(a_0)-
R(a_1)$.
\end{cor}

\section{The orientation of $\Sigma^1$}

Now let us recall a special case of Lemma 6.1 of \cite{ESz}. Let
$F \colon W^4 \to \Real^5$ be a generic map of a compact
orientable manifold. Then the singularity set $\Sigma(F)$ of $F$
is a 2-dimensional submanifold of $W^4$ which is not necessarily
orientable.

\begin{lem} \label{lem:1}
The line bundles $\det(T\Sigma(F))$ and $\ker(dF)$ over
$\Sigma(F)$ are isomorphic.
\end{lem}

\begin{defn}
Let $\pi$ denote the projection of $\Real^{m+1}$ onto $\Real^m$. A
map $f \colon N^n \to \Real^m$ is called \emph{prim} if there
exists an immersion $f' \colon N^n \looparrowright \Real^{m+1}$
such that $\pi \circ f' = f$.
\end{defn}

\begin{cor} \label{cor:1}
If $F \colon W^4 \to \Real^5$ is a generic prim map then
$\Sigma(F) \subset W^4$ is an orientable surface.
\end{cor}

\begin{proof}
Let $s$ denote the sixth coordinate function of $F'$, i.e., $F' =
(F, s)$. Since $F'$ is non-singular, the function $s$ is
non-degenerate along $\ker(dF)$. Thus we can orient $\ker(dF)$ so
that the derivative of $s$ in the positive direction of $\ker(dF)$
is positive. But the orientability of $\ker(dF)$ implies the
orientability of $\Sigma(F)$ by Lemma \ref{lem:1}.
\end{proof}

The following definition, motivated by Corollary \ref{cor:1},
gives an explicit isomorphism $\Psi$ between $\ker(dF)$ and
$\det(T\Sigma(F))$.

\begin{defn} \label{defn:4}
Let $W^4$ be a compact oriented manifold with possibly non-empty
boundary and let $F \colon W^4 \to \Real^5$ be a generic map. For
$p \in \Sigma(F)$ choose a small neighborhood $U_p \subset W^4$ of
$p$ in which $\ker(dF)$ is orientable. Put $F_p = F|U_p$ and
choose an orientation $o_p$ of $\ker(dF_p)$. Then there exists a
smooth function $s \colon U_p \to \Real$ such that the derivative
of $s$ in the direction of $o_p$ is positive. (First construct $s$
along $\Sigma(F_p)$ near $\Sigma^{1,1}(F_p)$ then extend it to a
tubular neighborhood of $\Sigma(F_p)$.) The map $F_p' = (F_p, s)
\colon U_p \looparrowright \Real^6$ is an immersion. If $U_p$ is
chosen sufficiently small then we can even suppose that $F_p'$ is
an embedding. Denote by $e_6$ the sixth coordinate direction in
$\Real^6$ and let $\nu_6 \colon U_p \to T\Real^6$ denote the
vector field along $F_p'$ defined by the formula $\nu_6(x) = e_6
\in T_{F_p'(x)}\Real^6$ for $x \in U_p$. Projecting $\nu_6$ into
the normal bundle of $F_p'$ we obtain a normal field $\mu_6$ along
$F_p'$ that vanishes exactly at the points of $\Sigma(F_p)$.
Perturb $F_p'$ in the direction of $\mu_6$ to obtain an embedding
$F_p''$. Then orient $\Sigma(F_p)$ as the intersection of $F_p'$
and $F_p''$ in $\Real^6$. Here $\Real^6$ is considered with its
standard orientation. This orientation of $\Sigma(F_p)$ does not
depend on the choice of the function $s$, since if $s_1$ and $s_2$
are two such functions then for $0 \le t \le 1$ the convex
combination $(1-t)s_1 + ts_2$ also satisfies the conditions for
$s$.

If we reverse the orientation of $\ker(dF_p)$, i.e. if we orient
it by $-o_p$, then we can choose $-s$ instead of $s$. Thus we
obtain the embedding $(F_p, -s)$, which is the reflection of
$F_p'$ in the hyperplane $\Real^5$. Denote this reflection by $R
\colon \Real^6 \to \Real^6$ (i.e., $R(x_1, \dots, x_5, x_6) =
(x_1, \dots, x_5, -x_6)$). Then $(F_p,-s) = R \circ F_p'$. The
vector field $dR \circ \nu_6$ along $R \circ F_p'$ points in the
direction $-e_6$ and $R \circ F_p''$ is the perturbation of $R
\circ F_p'$ in the direction of $dR \circ \mu_6$. But in this case
we should perturb $R \circ F_p'$ in the direction of $-(dR \circ
\mu_6)$. We obtain the same orientation if we look at the
intersection $(R \circ F_p'') \cap (R \circ F_p')$ instead. Since
the intersection is 2-dimensional and $U_p$ is 4-dimensional we
get that $(R \circ F_p'') \cap (R \circ F_p') = (R \circ F_p')
\cap (R \circ F_p'')$ in the oriented sense. The orientations of
$\Sigma(F_p)$ defined by the intersections $ F_p' \cap F_p''$ and
$(R \circ F_p') \cap (R \circ F_p'')$ are opposite. This can be
seen from the following argument: For $0 \le t \le 1$ denote by
$R_t$ the rotation of the hyperplane $\Real^6$ in $\Real^7$ around
$\Real^5$ by the angle $\pi t$. The orientation of $(R_t \circ
F_p')\cap (R_t \circ F_p'')$ in $R_t(\Real^6)$ changes
continuously as $t$ goes from $0$ to $1$. The orientations of the
hyperplanes $R_1(\Real^6)$ and $\Real^6$ are opposite, thus the
reflection $R$ changes the orientation of the intersection $F_p'
\cap F_p''$.

So we have defined an isomorphism $\Psi_p$ between $\ker(dF_p)$
and $\det(T\Sigma(F_p))$ for every $p \in W^4$ in a compatible way
(i.e., $\Psi_p|(U_p \cap U_q) = \Psi_q|(U_p \cap U_q)$ for $p,q
\in W^4$). These local isomorphisms define a global isomorphism
$\Psi$ between $\ker(dF)$ and $\det(T\Sigma(F))$.
\end{defn}

\section{The invariant}

In this section we will give a geometric formula for the
3-dimensional obstruction to the existence of a regular homotopy
between two immersions of $M^3$ into $\mathbb{R}^5$. This
generalizes the results of \cite{Takase} to immersions with
non-trivial normal bundle.

\begin{defn} \label{defn:2}
Let $W^4$ be a compact oriented manifold with boundary $M^3$ and
$F \colon W^4 \to \Real^5$ a generic map such that $f = F|M^3$ is
an immersion. Recall that $\Sigma(F)$ denotes the set of singular
points of $F$. Let us denote by $C(F) \subset M^3$ the
1-dimensional submanifold $\partial \Sigma(F)$. Choose a
trivialization $\tau$ of $\ker(dF)|C(F)$ so that it points into
the interior of $W^4$. This is possible since $f$ is non-singular
(and so $\ker(dF)$ never lies in $TM^3$). Then $\tau$ is normal to
$\Sigma(F)$ because $F$ is generic and thus $\Sigma^{1,1}(F) \cap
C(F) = \emptyset$. So if we project $\tau$ into $TM$ along
$\Sigma(F)$ we obtain a nowhere vanishing normal field $\nu(F)$ in
$\nu(C(F) \subset M^3)$.

Let $U$ denote a small collar neighborhood of $C(F)$ in
$\Sigma(F)$. Then clearly $U$ is orientable. Using the isomorphism
$\Psi$ of Definition \ref{defn:4} the trivialization $\tau$ of
$\ker(dF)$ induces an orientation of $U$. Thus $C(F) \subset
\partial U$ is also oriented. So we have assigned a pair $(C(F),
\nu(F))$ to $F$ as in Notation \ref{note:1}. Let $r(F) = r(C(F),
\nu(F))$.
\end{defn}

\begin{note}
For $\chi \in H^2(M^3; \mathbb{Z})$ let us denote by $\imm(M^3,
\Real^5)_{\chi}$ the space of immersions with normal Euler class
$\chi$.
\end{note}

Fix a cohomology class $\chi \in H^2(M^3; \mathbb{Z})$. Our aim is
to define an invariant $i \colon \pi_0\left(\imm(M^3,
\Real^5)_{\chi}\right) \to \mathbb{Z}_{2d(\chi)}$.

\begin{prop}
Let $f \in \imm(M^3, \Real^5)_{\chi}$ and let $F \colon W^4 \to
\Real^5$ be a generic map such that $\partial F = f$. Then $[C(F)]
= D\chi$.
\end{prop}

\begin{proof}
Let $\kappa$ denote an inner normal field of $M^3$ in $W^4$ that
extends $\tau$ (see Definition \ref{defn:2}). Then $dF \circ
\kappa$ is a vector field along $f$ that is tangent to $f$ exactly
at the points of $C(F) =
\partial \Sigma(F)$. (If $p \in C(F)$ then the rank of $(dF)_p$ is
$3$, moreover $dF|(T_pM^3) = df$ is non-degenerate. Thus
$dF(\kappa_p) \in df(T_pM^3)$.) So if we project $dF \circ \kappa$
into the normal bundle of $f$ we obtain a normal field of $f$ that
vanishes along $C(F)$. To see that $C(F)$ represents the normal
Euler class of $f$, we have to know that it is oriented suitably.

Using the notations of Definition \ref{defn:2} we choose a
function $s \colon W^4 \to \Real$ such that the derivative of $s$
in the direction of $\kappa$ (and thus $\tau$) is positive and $
s|M^3 \equiv 0 $. Then there exists a collar neighborhood $V$ of
$M^3$ in $W^4$ such that $F' = (F,s)|V$ is an immersion. Denote by
$\nu_{F'}$ the normal bundle of $F'$ in $\Real^6$ and by $\nu_f$
the normal bundle of $f$ in $\Real^5$. Then $\nu_{F'}|M^3 = \nu_f$
as oriented bundles, since $s$ is increasing along $\kappa$ (here
$\Real^5$ and $\Real^6$ are considered with their standard
orientations). By Definition $\ref{defn:4}$ the surface of
singular points $U = \Sigma(F|V)$ is oriented as the
self-intersection of $F'$ in $\Real^6$, or more precisely, as the
intersection of the zero section and a generic section of
$\nu_{F'}$. Moreover, $C(F)$ is oriented as the boundary of $U$.
Thus $C(F)$ is the self-intersection of the zero section of
$\nu_{F'}|M^3 = \nu_f$, so it is dual to the Euler class $e(\nu_f)
= \chi$. (Here we used the naturality of the Euler class.)
\end{proof}

\begin{defn} \label{defn:3}
Let $f \in \imm(M^3, \Real^5)_{\chi}$ and let $F \colon W^4 \to
\Real^5$ be generic such that $\partial F = f$. Denote the
algebraic number of cusps of $F$ by $\#\Sigma^{1,1}(F)$ (for the
definition see \cite{ESz}). Then let $j(f) = 3\sigma(W^4) - 3
\alpha(M^3) + \#\Sigma^{1,1}(F) + r(F) \in \mathbb{Z}_{2d(\chi)}$.
\end{defn}

Note that if $\chi = 0$ and $F$ is an immersion near $\partial
W^4$ then $r(F) = 0$. Thus in this case $j(f)$ agrees with the
double of the invariant introduced in \cite{Takase} (see
Definition \ref{defn:6}).

\begin{thm} \label{thm:3}
j(f) is well defined, i.e., it does not depend on the choice of
the generic map $F$. Moreover, if $f_0$ and $f_1$ are regularly
homotopic then $j(f_0) = j(f_1)$.
\end{thm}

\begin{proof}
For $i \in \{0,1\}$ let $F_i \colon W^4_i \to \Real^5$ be a
generic map such that $\partial F_i = f_i$. Choose a regular
homotopy $\{h_t: 0 \le t \le 1\}$ connecting $f_0$ and $f_1$. This
defines an immersion $H \colon M^3 \times I \to \Real^5 \times I$
by the formula $H(x,t) = (h_t(x),t)$. Also choose a closed collar
neighborhood $U_i$ of $M^3$ in $W^4_i$ and a diffeomorphism $d_i
\colon U_i \to M^3 \times [0,\varepsilon]$ for $i = 0,1$. Let $p
\colon M^3 \times [0, \varepsilon] \to [0, \varepsilon ]$ denote
the projection onto the second factor. If $\varepsilon$ (i.e.,
$U_i$) is sufficiently small then $p \circ d_i$ is non-degenerate
along $\ker(dF_i)$ for $i =0,1$ since $\ker(dF_i)$ never lies in
$TM^3$. Let $s_i$ be an arbitrary smooth extension of $p \circ
d_i$ over $W^4_i$. Now let $F_0' = (F_0, -s_0) \colon -W^4_0 \to
\Real^6$ and $F_1' = (F_1, s_1 + 1) \colon W^4_1 \to \Real^6$.
Then $F_i'$ is an immersion on $U_i$. Notice that $H|(M^3 \times
\{0\}) = F_0'| (\partial W^4_0)$ and $H|(M^3 \times \{1\}) = F_1'|
(\partial W^4_1)$.

Denote by $\kappa_i$ the inner normal field of $W^4_i$ along $M^3
= \partial W^4_i$ and by $v_6$ the sixth coordinate direction in
$\Real^6$. Then the inner product $\langle dF_0'(\kappa_0), v_6
\rangle < 0$ and $\langle dF_1'(\kappa_1), v_6 \rangle > 0$.
Furthermore, if $\lambda_i$ denotes the inner normal field of $M^3
\times I$ along $M^3 \times \{i\}$ for $i = 0,1$ then $\langle
dH(\lambda_0), v_6 \rangle > 0$ and $\langle dH(\lambda_1), v_6
\rangle < 0$. So $dH(\lambda_i)$ is homotopic to $-dF_i(\kappa_i)$
in the space of vector fields normal to $H|(M^3 \times \{i\})$.
Using Smale's lemma there exists a regular homotopy of $H$ fixed
on the boundary $M^3 \times \{0,1\}$ that induces the above
homotopy of normal fields. Denote by $H'$ the result of this
regular homotopy of $H$. Then $F_0'$, $H'$ and $F_1'$ fit together
to a smooth map $F'$ of $W^4 = -W^4_0 \cup (M^3 \times I) \cup
W^4_1$ into $\Real^6$ that is an immersion on $M^3 \times I$. Let
$\pi \colon \Real^6 \to \Real^5$ denote the projection map. Then
by a small perturbation of $H'$ we can achieve that $F = \pi \circ
F'$ is generic.

Since $G = F | (M^3 \times I) = \pi \circ H'$ is prim, the
singular surface $\Sigma(G)$ is oriented and a trivialization
$\tau$ of $\ker(dG)$ is given. If we project $\tau$ into
$\nu(\Sigma(G) \subset M^3 \times I)$ we obtain a normal field
$\nu$ along $\Sigma(G)$ that vanishes exactly at the cusps of $G$,
i.e., where $\tau$ is tangent to $\Sigma(G)$. So
$\#\Sigma^{1,1}(G)$ is equal to the algebraic number of zeroes of
$\nu$, which in turn is congruent to $\text{rd} ((C(F_0),
\nu(F_0)),(C(F_1), \nu(F_1))) = r(F_0)-r(F_1)$ modulo $2d(\chi)$
by Definition \ref{defn:1}.

Now using the result of Sz\H{u}cs \cite{Szucs} that $3\sigma(W^4)
+ \#\Sigma^{1,1}(F) = 0$ we get that
\begin{multline} \label{eqn:1}
-\left(3\sigma(W_0^4) + \#\Sigma^{1,1}(F_0) + r(F_0)\right) +
\left(3\sigma(W_1^4) + \#\Sigma^{1,1}(F_1) + r(F_1)\right) = \\
3\sigma(W^4) + \#\Sigma^{1,1}(-F_0 \cup G \cup F_1) = 0.
\end{multline}
In the special case $f_0 = f_1 = f$ this implies that $j(f)$ is
well defined, and for $f_0$ and $f_1$ arbitrary (but regularly
homotopic) we get that $j$ is a regular homotopy invariant.
\end{proof}

\begin{prop} \label{prop:4}
For any immersion $f \colon M^3 \looparrowright \Real^5$ the
invariant $j(f)$ is always an even element of
$\mathbb{Z}_{2d(\chi)}$, i.e., it is mapped to $0$ by the
epimorphism $\mathbb{Z}_{2d(\chi)} \to \mathbb{Z}_2$.
\end{prop}

\begin{proof}
Choose an immersion $f_1 \in \imm(M^3, \Real^5)_0$ and denote $f$
by $f_0$. Since in \cite{Takase} it is proved that $j(f_1)$ is
always even ($j(f_1) = 2i(f_1)$ for $i$ as in Definition
\ref{defn:6}) it is sufficient to prove that $j(f_0) \equiv j(f_1)
\mod 2$. Choose a singular Seifert surface $F_i \colon W_i \to
\Real^5$ for $f_i$ ($i = 0,1$) and let $G \colon M^3 \times I \to
\Real^5$ be a generic map such that $-F_0 \cup G \cup F_1$ is a
smooth map on $-W_1 \cup (M^3 \times I) \cup W_2$. Then by
equation \ref{eqn:1} above it is sufficient to prove that
$r(F_0)-r(F_1) \equiv \#\Sigma^{1,1}(G) \mod 2$. Since $f_1$ has
trivial normal bundle we may choose $F_1$ to be an immersion in a
neighborhood of $\partial W^4_1$. So $G$ is an immersion in a
neighborhood of $M^3 \times \{1\}$, moreover $r(F_1) = 0$. The
difference between the present situation and the proof of Theorem
\ref{thm:3} is that now $\ker(dG)$ might be non-orientable. Using
Definition \ref{defn:7} we get that $r_2(F_0) - r_2(F_1) \equiv
r(F_0)-r(F_1) \mod 2$. Let $\nu$ denote a generic normal field
along $\Sigma^1(G)$ that extends both $\nu(F_0)$ and $\nu(F_1)$.
By definition $r_2(F_0)-r_2(F_1)$ equals the mod $2$ number of
zeroes of $\nu$. Thus we only have to prove that
$\abs{\nu^{-1}(0)} \equiv \#\Sigma^{1,1}(G) \mod 2$.

From now on we will denote $\Sigma^1(G)$ by $K$ and the line
bundle $\ker(dG) < T(M^3 \times I)|K$  by $l$. Then $l$ is tangent
to $K$ exactly at the points of $\Sigma^{1,1}(G)$. For $\eps >0$
sufficiently small let $\widetilde{K}$ denote the sphere bundle
$S_{\eps}l$. If $\eps$ is sufficiently small then the exponential
map of $M^3 \times I$ defines an immersion $s \colon \widetilde{K}
\to M^3 \times I$ so that the double points of $s$ correspond
exactly to the points of $\Sigma^{1,1}(G)$. So we have to prove
that $\abs{D_2(s)} \equiv \abs{\nu^{-1}(0)} \mod 2$. By Lemma
\ref{lem:1} the surface $\widetilde{K}$ is the orientation double
cover of $K$, in particular $\widetilde{K}$ is oriented and a sign
can be given to each double point of $s$ (here we also use that
$\dim(\widetilde{K})$ is even). The sign of a double point of $s$
is the opposite of the sign of the corresponding cusp of $G$
(since the sign of a cusp is defined as the self intersection of
$K$). Thus $\#D_2(s) = -\#\Sigma^{1,1}(G)$. Let $p \colon
\widetilde{K} \to K$ denote the covering map. Then $p^*\nu_K
\approx \nu_s$, thus $p^*\nu$ defines a section $\widetilde{\nu}$
of $\nu_s$. From the construction of $\widetilde{\nu}$ it is clear
that $\#\widetilde{\nu}^{-1}(0) = 2 \#\nu^{-1}(0)$.

If we perturb $s$ in the direction of $\widetilde{\nu}$ we get a
self intersection point of $s$ for each element of
$\widetilde{\nu}^{-1}(0)$ and two self intersection points for
each double point of $s$. Thus $$s \cap (s + \eps \widetilde{\nu})
=  \#\widetilde{\nu}^{-1}(0) + 2D_2(s) = 2(\#\nu^{-1}(0) -
\#\Sigma^{1,1}(G)).$$ So we only have to show that the left hand
side is divisible by $4$.

From now on we will work in a fixed tubular neighborhood $T$ of
$C(F_0) \subset M^3$. Note that $\partial (s, \widetilde{\nu}) =
(C(F_0)+ \eps \nu(F_0), \nu(F_0)) \cup (-C(F_0)-\eps \nu(F_0),
\nu(F_0)) \subset T \times \{0\}$. Let us denote by $C$ the
one-dimensional submanifold $(C(F_0)+ \eps \nu(F_0))\cup(-C(F_0)-
\eps \nu(F_0))$ of $M^3$. Then $s \cap (s + \eps\widetilde{\nu}) =
r(C, \nu(F_0)) \in \mathbb{Z}$ since $C \sim C(F_0)-C(F_0)$ is
null homologous in $M^3$. We define an embedding $e \colon C(F_0)
\times [-\eps, \eps] \to T $ by the formula $e(x,t) = x + t
\nu(F_0)$. Then $E=\im(e)$ is a 2-dimensional oriented submanifold
of $T$ with boundary $C$. Thus $r(C, \nu(F_0)) = E \cap
(C+\nu(F_0))$ where the right hand side is considered to be a
generic intersection (each fiber of $E$ is parallel to
$\nu(F_0)$). Let $n$ be a small non-zero vector field along
$C(F_0)$ orthogonal to $\nu(F_0)$. Then $E \cap (C + \nu(F_0)+n) =
\emptyset$ (this can be verified by inspecting each fiber of $T$),
thus $r(C, \nu(F_0)) = 0$. So we get that $\#\nu^{-1}(0) =
\#\Sigma^{1,1}(G)$, not just a mod 2 congruence.
\end{proof}

\begin{rem}
A small improvement on the proof of Proposition \ref{prop:4}
yields an interesting result: Let $G \colon M^3 \times [0,1] \to
\Real^5 $ be a generic map connecting the immersions $f_0$ and
$f_1$. Let $K$ denote the singular surface of $G$ and let $\nu_i$
be a trivialization of $\ker(dG)|(M^3 \times \{i\})$ for $i =
0,1$. Then $\#\Sigma^{1,1}(G)$ is equal to the relative twisted
normal Euler class $e(\nu_K; \nu_0, \nu_1)$.

Note that if $\partial K = C_0 \cup C_1$ then by definition
$$\text{rd}_2((C_0, \nu_0),(C_1,\nu_1)) \equiv e(\nu_K; \nu_0,\nu_1) \mod
2.$$ If in particular $K$ is an oriented cobordism between $C_0$
and $C_1$ then $$\text{rd}((C_0, \nu_0),(C_1,\nu_1)) = e(\nu_K;
\nu_0,\nu_1) = \#\Sigma^{1,1}(G);$$ here $C_i$ is oriented by
$\nu_i$ using the isomorphism $\Psi$ (see Definition
\ref{defn:4}).
\end{rem}

Thus $j$ may take only $d(\chi)$ different values if $d(\chi) >
0$. (Since $j$ is additive if a connected sum is taken with an
immersion of a sphere (Lemma \ref{lem:2}) and $j(g)$ can be any
even number for $g \colon S^3 \looparrowright \Real^5$ it follows
that $j$ is an epimorphism onto $2\mathbb{Z}_{2d(\chi)}$.) But
Theorem \ref{thm:1} implies that there are exactly $2d(\chi)$
regular homotopy classes with normal Euler class $\chi$. So $j$
describes the regular homotopy class of $f$ only up to a $2:1$
ambiguity. To resolve this problem we will lift the invariant $j
\in \mathbb{Z}_{2d(\chi)}$ to an invariant $I \in
\mathbb{Z}_{4d(\chi)}$. It follows from Proposition \ref{prop:4}
that $I$ is always an even element, thus it defines an invariant
$i \in \mathbb{Z}_{2d(\chi)}$ by the embedding
$\mathbb{Z}_{2d(\chi)} \hookrightarrow \mathbb{Z}_{4d(\chi)}$.

\begin{note}
Let $f \in \imm(M^3, \Real^5)_{\chi}$ and let $F \colon W^4 \to
\Real^5$ be a singular Seifert surface for $f$. If $\kappa$
denotes the inner normal field along $\partial W^4$ then let
$\bar{w}(F) \in \Gamma(\nu_f)$ be the projection of $dF(\kappa)$
into $\nu_f$. If $\bar{t} \in \Gamma(\eps^1_M)$ denotes a
trivialization of $\eps^1_M$ then we can consider $\bar{t}$ and
$\bar{w}(F)$ to be sections of $\nu_f \oplus \eps^1_M$. Let
$\bar{v}(F) = \bar{w}(F) + \bar{t} \in \Gamma(\nu_f \oplus
\eps^1_M)$.
\end{note}

From now on we will fix a spin structure $s_M \in
\text{Spin}(M^3)$. If we consider $\Real^5$ with its unique spin
structure then for every immersion $f \colon M^3 \looparrowright
\Real^5$ a spin structure $s(f)$ is induced on $\nu_f$ by $s_M$.
Then $s(f)$ is equivalent to a trivialization $\tau(f) \colon
\eps^3_M|\text{sk}_2(M) \to \nu_f \oplus
\eps^1_M|\text{sk}_2(M^3)$ up to homotopy. Since $\pi_2(SO(3)) =
0$ the trivialization $\tau(f)$ extends to an isomorphism $\tau(f)
\colon \eps^3_M \to \nu_f \oplus \eps^1_M$, but this extension is
not unique because $\pi_3(SO(3)) \neq 0$.

\begin{defn}
Using the above notations let $t(f),v(F) \in \Gamma(\eps^3_M)$ be
defined by the formulas $t(f) = \tau(f)^{-1} \circ \bar{t}$ and
$v(F) = \tau(f)^{-1} \circ \bar{v}(F)$. Denote by $a(F)$ the
quadruple $(C(F), \nu(F), t(f), v(F)) \in N(M^3, \chi)$. Then
define $R(F) \in \mathbb{Z}_{4d(\chi)}$ to be $R(a(F))$. Since the
homotopy class of the map $(t(f),v(F))|\text{sk}_2(M^3) \colon
\text{sk}_2(M^3) \to V_2(\Real^3)$ is independent of the choice of
the extension of $\tau(f)$ to $\eps^3_M$ Proposition \ref{prop:3}
implies that $R(F)$ is also independent of $\tau(f)$ and depends
only on $s_M$.
\end{defn}

\begin{rem} \label{rem:2}
Proposition \ref{prop:5} implies that $r(F) \equiv R(F) \mod
2d(\chi)$.
\end{rem}

Now we can finally define a complete regular homotopy invariant.

\begin{defn}
For $f \in \imm(M^3, \Real^5)_{\chi}$ and a singular Seifert
surface $F$ let $I(f) \in \mathbb{Z}_{4d(\chi)}$ be defined as
$3\sigma(W^4) - 3 \alpha(M^3) + \#\Sigma^{1,1}(F) + R(F)$. (Recall
that we have fixed a spin structure $s_M$ on $M^3$ for the
definition of $R(F)$.) Remark \ref{rem:2} above implies that $j(f)
\equiv I(F) \mod 2d(\chi)$. Thus by Proposition \ref{prop:4} we
get that $I(F)$ is always an even element of
$\mathbb{Z}_{4d(\chi)}$. Let us denote by $\frac12$ the
isomorphism from $2\mathbb{Z}_{4d(\chi)}$ to
$\mathbb{Z}_{2d(\chi)}$. Then let $i(F) = \frac12 I(F)$.
\end{defn}

Clearly $j(f) = 2i(f)$ for every $f \in \imm(M^3,
\Real^5)_{\chi}$.

\begin{thm}
I(f) is well defined, i.e., it does not depend on the choice of
the generic map $F$. Moreover, if $f_0$ and $f_1$ are regularly
homotopic then $I(f_0) = I(f_1)$.
\end{thm}

\begin{proof}
Using the notations of the proof of Theorem \ref{thm:3} we only
have to show that the surface $K = \Sigma(G) \subset M^3 \times I$
satisfies Definition \ref{defn:5}. I.e., there exist generic
sections $t$ and $v$ of $\eps^3_{M \times I}$ that extend $t(f_i)$
and $v(F_i)$ for $i = 0,1$ and are linearly dependent exactly over
$K$. The regular homotopy between $f_0$ and $f_1$ defines the
immersion $H \colon M^3 \times I \looparrowright \Real^5 \times
I$. For $i \in \{0,1\}$ there is a canonic isomorphism $\varphi_i
\colon \nu_H|(M^3 \times \{i\}) \to \nu_{f_i}$. Let $\bar{w} \in
\Gamma(\nu_H)$ denote the projection of the sixth coordinate
vector $v_6 \in \Real^6$ into $\nu_H$. Then $\varphi_i \circ
(\bar{w}|M^3 \times \{i\}) = \bar{w}(F_i)$. Moreover, $K =
\bar{w}^{-1}(0)$ and the orientation of $K$ is defined as the self
intersection of $H$ if perturbed in the direction of $\bar{w}$.
Define $\bar{t}$ to be a trivialization of the $\eps^1_{M \times
I}$ component of the bundle $\nu_H \oplus \eps^1_{M \times I}$ and
let $\bar{v} = \bar{w} + \bar{t}$. Then $\bar{t}$ and $\bar{v}$
are linearly dependent exactly at the points of $K$. Note that
$t(f_i) = \tau(f_i)^{-1} \circ (\varphi_i \oplus
\text{id}_{\eps^1})(\bar{t}|M^3 \times \{i\})$ and $v(F_i) =
\tau(f_i)^{-1} \circ (\varphi_i \oplus
\text{id}_{\eps^1})(\bar{v}|M^3 \times \{i\})$. Thus we only have
to define a trivialization $\tau \colon \eps^3_{M \times I} \to
\nu_H \oplus \eps^1_{M \times  I}$ such that
\begin{equation} \label{eqn:2}
\tau|(M^3 \times \{i\}) = \varphi_i^{-1} \circ \tau(f_i)
\,\,\,\,\text{for}\,\,\,\, i = 0,1.
\end{equation}

The spin structure $s_M \in \text{Spin}(M^3)$ and the unique spin
structure on $I$ define a spin structure on $M^3 \times I$.
Together with the unique spin structure of $\Real^6$ we get a spin
structure $s_H$ on $\nu_H$. When $s_H$ is restricted to $M^3
\times \{i\}$ we get back the spin structure $s_M$. Thus $s_H$
defines a trivialization $$\tau_H \colon \eps^3_{M \times
I}|\text{sk}_2(M^3 \times I) \to (\nu_H \oplus \eps^1_{M \times
I})|\text{sk}_2(M^3 \times I)$$ satisfying equation \ref{eqn:2}
over the 2-skeleton of $M^3 \times \{0,1\}$. Note that the
trivialization $\tau(f_i)$ is only well defined over
$\text{sk}_2(M^3)$ and that we can choose an arbitrary extension
over $M^3$ in order to define the rotation difference. Thus we
only have to extend $\tau_H$ to a trivialization $\tau$ of $\nu_H
\oplus \eps^1_{M \times I}$ and then define $\tau(f_i)$ by formula
\ref{eqn:2}.

First we extend $\tau_H$ to $\text{sk}_3(M^3 \times I) \setminus
\text{sk}_3(M^3 \times \{1\})$. This is possible since the
obstruction to extending the trivialization over a 3-simplex from
its boundary lies in $\pi_2(SO(3)) = 0$. If $\sigma^3$ is a
3-simplex of $M^3$ then we can extend $\tau_H$ to $\sigma^3 \times
I$ since it is given only on $\partial (\sigma^3 \times I)
\setminus (\sigma^3 \times \{1\})$. Thus we have obtained the
required extension $\tau$ of $\tau_H$.
\end{proof}

\section{Connected sums and completeness of the invariant $i$}

\begin{lem} \label{lem:2}
If $f \in \imm(M^3, \Real^5)_{\chi}$ and $g \in \imm(S^3,
\Real^5)$ then $c(f \# g) = c(f)$, in particular $e(\nu_{f \# g})
= e(\nu_f)$. Moreover
$$i(f \# g) = i(f) + \left(i(g) \mod 2d(\chi)\right) \in
\mathbb{Z}_{2d(\chi)}.
$$
\end{lem}

\begin{proof}
Since $c(f)$ describes the regular homotopy class of
$f|(M^3_{\circ})$ it is trivial that $c(f \# g) = c(f)$. Let $F$
be a singular Seifert surface of $f$ and $G$ of $g$ such that $G$
is an immersion near the boundary. Then the result follows by
inspecting the boundary connected sum $F \natural G$ and the fact
that $C(G) = \emptyset$.
\end{proof}

\begin{thm}
Suppose that the immersions $f_0, f_1 \in \imm(M^3, \Real^5)$ are
regularly homotopic on $M^3 \setminus D$, where $D \subset M^3$ is
diffeomorphic to the closed disc $D^3$ (i.e., $c(f_0) = c(f_1)$).
Then $i(f_0) = i(f_1)$ implies that $f_0$ is regularly homotopic
to $f_1$.
\end{thm}

\begin{proof}
The proof consists of two cases according to the value of
$d(\chi)$.

If $d(\chi) > 0$ then $i$ takes values in $\mathbb{Z}_{2d(\chi)}$
which is a finite group. Theorem \ref{thm:1} implies that there
are exactly $2d(\chi)$ regular homotopy classes with a fixed Wu
invariant $c$. Thus we only have to show that the invariant $i$
restricted to immersions with Wu invariant $c$ is an epimorphism
onto $\mathbb{Z}_{2d(\chi)}$. For this end choose an immersion $f
\in \imm(M^3, \Real^5)$ such that $c(f) = c$. In \cite{ESz} it is
shown that $i \colon \imm[S^3, \Real^5] \to \mathbb{Z}$ is a
bijection. Thus Lemma \ref{lem:2} implies that $c(f \# g) = c(f) =
c$ for every $g \in \imm(S^3, \Real^5)$, moreover $i \colon \{\, f
\#g \colon g \in \imm(S^3, \Real^5) \,\} \to
\mathbb{Z}_{2d(\chi)}$ is surjective.

If $d(\chi) = 0$ then $i$ maps into $\mathbb{Z}$. Using Smale's
lemma we can suppose that $f_0|(M^3 \setminus D) = f_1|(M^3
\setminus D)$. The normal bundles of $f_0|D$ and $f_1|D$ in
$\Real^5$ are trivial, choose a trivialization for both of them.
Let $\tau_0$ be a non-zero normal field along $f_0|D$. Then
$\tau_0| \partial D$ considered in the trivialization of the
normal bundle of $f_1|D$ is a map $(\tau_0|
\partial D) \colon
\partial D \to S^1$. Since $\partial D$ is homeomorphic to $S^2$
and $\pi_2(S^1) = 0$ the normal field $\tau_0|\partial D$ can be
extended to a normal field $\tau_1$ of $f_1|D$. Thus $\tau_i$ is a
normal field of $f_i|D$ for $i = 0,1$ and $\tau_0|\partial D =
\tau_1 | \partial D$.

Next choose an oriented compact manifold $W^4_0$ with boundary
$M^3$. We push $D$ into the interior of $W^4_0$ fixing the
boundary $\partial D$ to obtain a 3-disc $D_1 \subset W^4_0$ so
that $\partial D = \partial D_1$ and $M^3_1 = (M^3 \setminus D)
\cup D_1$ is a smooth submanifold of $W^4_0$. If we throw out the
domain bounded by $D$ and $D_1$ in $W^4_0$ we obtain a
4-dimensional submanifold $W^4_1$ of $W^4_0$ with boundary
$M^3_1$. Clearly $W^4_0$ is diffeomorphic to $W^4_1$.

We can choose a generic map $F_0 \colon W^4_0 \to \Real^5$ with
the following three properties:

\begin{enumerate}
\item $F_0 | M^3 = f_0$ and $F_0 | M^3_1 = f_1$ (where $M^3_1$ is
identified with $M^3$ by a diffeomorphism keeping $M^3 \setminus
D$ fixed).
\item $F_0$ is an immersion in a neighborhood of $D$ and $D_1$.
\item If $\kappa_0$ denotes the inner normal field of $D$ in
$W^4_0$ and $\kappa_1$ denotes the inner normal field of $D_1$ in
$W^4_1$ then $dF_0 \circ \kappa_0 = \tau_0$ and $dF_1 \circ
\kappa_1 = \tau_1$.
\end{enumerate}

Let $F_1 = F_0 | W^4_1$. Then (2) implies that $C(F_0) = C(F_1)
\subset M^3 \setminus D$, moreover $\nu(F_0) = \nu(F_1)$. In
particular, the normal Euler class of $f_0$ and $f_1$ coincide.
Thus $R(F_0) = R(F_1)$. Since $\sigma(W^4_0) = \sigma(W^4_1)$, we
get that $$0 = i(f_0)-i(f_1) = \#\Sigma^{1,1}(F_0|(W^4_0 \setminus
W^4_1)).$$

Choose diffeomorphisms $d_0 \colon S^3_+ \to D$ and $d_1 \colon
S^3_+ \to D_1$, where $S^3_+$ denotes the northern hemisphere of
$S^3$. Then the immersion $F_0 \circ d_i$ can be extended to an
immersion $f_i' \colon S^3 \looparrowright \Real^5$ for $i = 0,1$
so that $f_0'|S^3_- = f_1'|S^3_-$. (This is possible since
$j^1(f_0)|\partial D = j^1(f_1)|\partial D.$) Now repeating the
same argument as above for $f_0'$ and $f_1'$, we obtain that
$$i(f_0')-i(f_1') =  \#\Sigma^{1,1}(F_0|(W^4_0 \setminus
W^4_1)).$$ (Note that $\tau_0$ and $\tau_1$ have a common
extension over $S^3_-$.) Thus $i(f_0') -i(f_1') = 0$, so using
\cite{ESz} we get that $f_0'$ and $f_1'$ are regularly homotopic.
But this implies that there exists a regular homotopy between
$f_0'$ and $f_1'$ that is fixed on $S^3_-$ (see \cite{Juhasz2},
Lemma 3.33). So $f_0|D$ and $f_1|D$ are regularly homotopic
keeping the 1-jets on the boundary fixed, which completes the
proof that $f_0$ and $f_1$ are regularly homotopic.
\end{proof}

\begin{cor}
The map $$(c,i) \colon \imm[M^3, \Real^5] \to \coprod_{c \in
H^2(M^3; \mathbb{Z})} \mathbb{Z}_{4d(c)}$$ is a bijection.
\end{cor}

We get more structure on the set of regular homotopy classes of
immersions of oriented 3-manifolds into $\Real^5$ if we endow it
with the connected sum operation. Let us introduce the notation
$$I(3,5) = \left\{\, [f] \colon [f] \in \imm[M^3,\Real^5] \,\,\text{for}\,\,
M^3 \,\,\text{oriented} \,\right\}.$$ Then $(I(3,5), \#)$ is a
semigroup whose structure is described in the following theorem.

\begin{thm}
Let $M^3_1$ and $M^3_2$ be oriented 3-manifolds. Then
\begin{equation} \label{eqn:3}
H^2(M^3_1 \# M^3_2; \mathbb{Z}) \approx H^2(M^3_1;
\mathbb{Z})\oplus H^2(M^3_2; \mathbb{Z}).
\end{equation}
If $f_i \in \imm(M^3_i, \Real^5)$ for $i = 1,2$ then
\begin{equation} \label{eqn:4}
c(f_1 \# f_2) = c(f_1) \oplus c(f_2) \in H^2(M^3_1 \# M^3_2;
\mathbb{Z}).
\end{equation}
Moreover, if $\chi_i$ denotes the normal euler class of $f_i$ and
$\chi$ the normal euler class of $f_1 \# f_2$ then $d(\chi) =
\gcd(d(\chi_1), d(\chi_2))$. Finally,
\begin{equation} \label{eqn:5}
i(f_1 \# f_2) = (i(f_1) \mod 2d(\chi)) + (i(f_2) \mod 2d(\chi))
\in \mathbb{Z}_{2d(\chi)},
\end{equation}
where $i(f_i) \in \mathbb{Z}_{2d(\chi_i)}$ for $i = 1,2$.
\end{thm}

\begin{proof}
Equation \ref{eqn:3} follows from the fact that $H^2(M^3_i;
\mathbb{Z}) \approx H^2(M^3_i \setminus D^3; \mathbb{Z})$ (see the
long exact sequence of the pair $(M^3_i, M^3_i \setminus D^3)$)
and the Mayer-Vietoris exact sequence for $M^3_1 \# M^3_2 = (M^3_1
\setminus D^3) \cup (M^3_2 \setminus D^3)$.

Equation \ref{eqn:4} can be seen from the description of $c(f_i)$
as the regular homotopy class of $f_i|\text{sk}_2(M^3)$. Since
$\chi = \chi_1 \oplus \chi_2$ the statement about $d(\chi)$ is
trivial.

Finally, equation \ref{eqn:5} is obtained by taking the boundary
connected sum $F_1 \natural F_2$ of singular Seifert surfaces
$F_1$ and $F_2$ for $f_1$, respectively $f_2$.
\end{proof}

\section{Immersions of $M^3$ into $\Real^6$ with a normal field}

Let $\imm_1(M^3, \Real^6)$ denote the space of immersions of $M^3$
into $\Real^6$ with a normal field $\nu$. Moreover, let
$\imm_1[M^3, \Real^6] = \pi_0(\imm_1(M^3, \Real^6))$ be the set of
regular homotopy classes of such immersions with normal fields. If
we fix a trivialization of $TM^3$ then  Hirsch's theorem
\cite{Hirsch} implies that the natural map $\imm_1(M^3, \Real^6)
\to C(M^3, V_4(\Real^6))$ is a weak homotopy equivalence.

For $f \in \imm(M^3, \Real^5)$ let $\iota(f) \in \imm_1(M^3,
\Real^6)$ be the immersion $f$ with the constant normal field
defined by the sixth coordinate vector in $\Real^6$. Thus $\iota$
is an embedding of $\imm(M^3, \Real^5)$ into $\imm_1(M^3,
\Real^6)$. As a special case of Hirsch's compression theorem we
have the following proposition.

\begin{prop}
$\iota_* \colon \imm[M^3, \Real^5] \to \imm_1[M^3, \Real^6]$ is a
bijection.
\end{prop}

\begin{proof}
The embedding $\Real^5 \hookrightarrow \Real^6$ induces an
embedding $V_3(\Real^5) \hookrightarrow V_4(\Real^6)$ and thus a
map $\psi \colon [M^3, V_3(\Real^5)] \to [M^3, V_4(\Real^6)]$ that
makes the following diagram commutative.
\[
\begin{CD}
\imm[M^3, \Real^5] @>\iota_*>>
\imm_1[M^3, \Real^6] \\
@VVV @VVV \\
[M^3, V_3(\Real^5)] @>\psi>> [M^3, V_4(\Real^6)].
\end{CD}
\]
By Hirsch's theorem the vertical arrows are bijections, thus it is
sufficient to prove that $\psi$ is also a bijection. To see this
consider the fibration $V_3(\Real^5) \to V_4(\Real^6) \to S^5$.
Then from the homotopy exact sequence of this fibration we get
that the homomorphism $\pi_i(V_3(\Real^5)) \to
\pi_i(V_4(\Real^6))$ is an isomorphism for $i \le 3$ and this
implies that $\psi$ is a bijection.
\end{proof}

The natural forgetful map $\varphi \colon \imm_1(M^3, \Real^6) \to
\imm(M^3, \Real^6)$ is a Serre fibration.

\begin{prop}
For any immersion $f \colon M^3 \looparrowright \Real^6$ the
normal bundle $\nu_f$ is trivial.
\end{prop}

\begin{proof}
Since $M^3$ is spin the normal bundle $\nu_f$ is also spin, thus
it is trivial over the 2-skeleton of $M^3$. Such a trivialization
can be extended to the 3-simplices of $M^3$ because
$\pi_2(SO(3))=0$.
\end{proof}

So $\varphi$ is surjective and the fiber of $\varphi$ is homotopy
equivalent to $\Gamma(\nu_f) = C(M^3, S^2)$. Thus the end of the
homotopy exact sequence of $\varphi$ looks like as follows:
\[\pi_1(\imm(M^3, \Real^6)) \longrightarrow [M^3, S^2] \longrightarrow
\imm_1[M^3, \Real^6] \xrightarrow{\,\,\varphi_*} \imm[M^3,
\Real^6] \longrightarrow 0.\] By Hirsch's theorem there is a
bijection $\imm[M^3, \Real^6] \approx [M^3, V_3(\Real^6)]$. Since
$V_3(\Real^6)$ is 2-connected and $\pi_3(V_3(\Real^6)) \approx
\mathbb{Z}_2$ we get from obstruction theory that $[M^3,
V_3(\Real^6)] \approx H^3(M^3; \mathbb{Z}_2) \approx
\mathbb{Z}_2$. It is well known that for $f \in \imm(M^3,
\Real^6)$ the regular homotopy class of $f$ is determined by the
number of its double points $D(f)$ modulo 2. This gives a
geometric interpretation of the map $\varphi_*$: for $(f, \nu) \in
\imm_1(M^3, \Real^6)$ the regular homotopy invariant $\varphi_*(f,
\nu)$ is equal to $D(f)$ modulo 2.

How can we determine the value of $\varphi_* \circ \iota_*(g)$ for
a generic $g \in \imm(M^3, \Real^5) ?$ This question was answered
in \cite{Szucs2}, let us recall that result now. The
self-intersection set $A(g)$ of $g$ is a closed 1-dimensional
submanifold of $M^3$ and $g(A(g))$ is also a closed 1-dimensional
submanifold of $\Real^5$. We say that a component $C$ of $g(A(g))$
is non-trivial if the double cover $g|g^{-1}(C) \colon g^{-1}(C)
\to C$ is non-trivial, i.e., if $g^{-1}(C)$ is connected. Let the
number of non-trivial components be denoted by $\delta(g)$. In
\cite{Szucs2} Sz\H{u}cs proved the following.

\begin{thm}
Suppose that $f \colon M^3 \to \Real^6$ is a generic immersion and
that $\pi \colon \Real^6 \to \Real^5$ is a projection such that $g
= \pi \circ f$ is also a generic immersion. Then $$D(f) \equiv
\delta(g) \mod 2.$$
\end{thm}

Note that the immersions $f$ and $g$ above are regularly homotopic
in $\Real^6$. Thus for any generic $g$ we have that $\varphi_*
\circ \iota_* = \delta$.

Now we are going to determine the group $\pi_1(\imm(M^3,
\Real^6))$. Using Hirsch's theorem we get that it is isomorphic to
$\pi_1(C(M^3, V_3(\Real^6))) = [SM^3, V_3(\Real^6)]$. Here $SM^3$
denotes the suspension of $M^3$ and is a 4-dimensional CW complex.
The space $V_3(\Real^6)$ is 2-connected and $\pi_3(V_3(\Real^6))
\approx \mathbb{Z}_2$. Moreover, $\pi_4(V_3(\Real^6)) \approx 0$.
This can be seen as follows: From the homotopy exact sequence of
the fibration $V_3(\Real^6) \to V_4(\Real^7) \to S^6$ we get that
$\pi_4(V_3(\Real^6)) \approx \pi_4(V_4(\Real^7))$. It was shown by
Paechter \cite{Paechter} that for $k \ge 4$ the isomorphism
$\pi_k(V_k(\Real^{2k-1})) \approx 0$ holds if $k \equiv 0 \mod 4$.
Thus obstruction theory yields that $[SM^3, V_3(\Real^6)] \approx
H^3(SM^3; \mathbb{Z}_2) \approx H^2(M^3; \mathbb{Z}_2)$.

Putting together the above results we obtain the following
theorem.

\begin{thm} \label{thm:4}
The following sequence is exact: $$H^2(M^3; \mathbb{Z}_2)
\longrightarrow [M^3, S^2] \longrightarrow \imm[M^3, \Real^5]
\xrightarrow{\,\,\delta\,\,} \mathbb{Z}_2 \longrightarrow 0.
$$
\end{thm}

\begin{rem}
If we fix a trivialization of $TM^3$ then non-zero vector fields
(or equivalently, oriented 2-plane fields) on $M^3$ correspond to
maps $M^3 \to S^2$. Thus the set of homotopy classes of oriented
2-plane fields on $M^3$ is equal to $[M^3, S^2]$ which was
determined in Remark \ref{rem:3}. A geometric classification of
such oriented 2-plane fields, avoiding the use of a trivialization
of $TM^3$, was carried out by Gompf (see \cite{Gompf}, section 4).
A complete set of homotopy invariants, similar to those introduced
in our present paper, were obtained in \cite{Gompf}. Gompf's
result and the regular homotopy classification of immersions of
$M^3$ into $\Real^5$ are related by Theorem \ref{thm:4}.
\end{rem}

\section*{Acknowledgements}

I would like to take this opportunity to thank Professor Andr\'as
Sz\H{u}cs for our long and helpful discussions and for reading
earlier versions of this paper.

% ----------------------------------------------------------------
\bibliographystyle{amsplain}
\bibliography{topology}

\providecommand{\bysame}{\leavevmode\hbox to3em{\hrulefill}\thinspace}
\providecommand{\MR}{\relax\ifhmode\unskip\space\fi MR }
\providecommand{\MRhref}[2]{%
  \href{http://www.ams.org/mathscinet-getitem?mr=#1}{#2}
}
\providecommand{\href}[2]{#2}
\begin{thebibliography}{10}

\bibitem{Li}
Li~Banghe, \emph{On classification of immersions of n-manifolds in
  (2n-1)-manifolds}, Comment. Math. Helvetici \textbf{57} (1982), 135--144.

\bibitem{Ekholm2}
T.~Ekholm, \emph{Regular homotopy and {V}assiliev invariants of generic
  immersions ${S}^k \to \mathbb{R}^{2k-1}$, $k \ge 4$}, J. Knoth Theory
  Ramifications \textbf{7} (1998), no.~8, 1041--1064.

\bibitem{ESz}
T.~Ekholm and A.~Sz\H{u}cs, \emph{Geometric formulas for {S}male invariants of
  codimension two immersions}, Topology \textbf{42} (2003), no.~1, 171--196.

\bibitem{Gompf}
R.~E. Gompf, \emph{Handlebody construction of {S}tein surfaces}, Ann. Math.
  \textbf{148} (1998), 619--693.

\bibitem{Hirsch}
M.~W. Hirsch, \emph{Immersions of manifolds}, Trans. Am. Math. Soc. \textbf{93}
  (1959), 242--276.

\bibitem{Hughes}
J.F. Hughes and P.M. Melvin, \emph{The {S}male invariant of a knot}, Comment.
  Math. Helv. \textbf{60} (1985), 615--627.

\bibitem{Juhasz2}
A.~Juh\'asz, \emph{Regular homotopy classes of singular maps}, Proc. London
  Math. Soc. (submitted).

\bibitem{Paechter}
G.~F. Paechter, \emph{The groups $\pi_r({V}_{n,m})$}, Quart. J. Math. Oxford
  Ser. \textbf{7} (1956), no.~2, 249--268.

\bibitem{Pontrjagin}
L.~Pontrjagin, \emph{A classification of mappings of the three-dimensional
  complex into the two dimensional sphere}, Matematicheskii Sbornik \textbf{9}
  (1941), no.~2, 331--363.

\bibitem{Takase}
O.~Saeki, A.~Sz\H{u}cs, and M.~Takase, \emph{Regular homotopy classes of
  immersions of 3-manifolds into 5-space}, Manuscripta Math. \textbf{108}
  (2002), 13--32.

\bibitem{Smale2}
S.~Smale, \emph{Classification of immersions of spheres in {E}uclidean space},
  Ann. Math. \textbf{69} (1959), 327--344.

\bibitem{Szucs2}
A.~Sz{\H{u}}cs, \emph{Note on double points of immersions}, Manuscripta Math.
  \textbf{76} (1992), 251--256.

\bibitem{Szucs}
\bysame, \emph{On the singularities of hyperplane projections of immersions},
  Bull. London Math. Soc. \textbf{32} (2000), 364--374.

\bibitem{Wu}
Wu{,} Wen-Ts{\"u}n, \emph{On the immersion of ${C}^{\infty}$-3-manifolds in a
  {E}uclidean space}, Sci. Sinica \textbf{13} (1964), 335--336.

\end{thebibliography}
\end{document}